\documentclass[reqno,11pt]{amsart}
\usepackage{amscd,amssymb}

\theoremstyle{plain}
\newtheorem{Thm}[subsection]{Theorem}
\newtheorem{Cor}[subsection]{Corollary}
\newtheorem{Lem}[subsection]{Lemma}
\newtheorem{Prop}[subsection]{Proposition}
\newtheorem{Conj}[subsection]{Conjecture}

\theoremstyle{definition}
\newtheorem{Def}[subsection]{Definition}

\theoremstyle{remark}

\newtheorem{Rem}[subsection]{Remark}

\numberwithin{equation}{section}
\renewcommand{\rm}{\normalshape}

\newif\ifShowLabels
\ShowLabelstrue
\newdimen\theight
\def\TeXref#1{%
    \leavevmode\vadjust{\setbox0=\hbox{{\tt
        \quad\quad  {\small \rm #1}}}%
    \theight=\ht0
    \advance\theight by \lineskip
    \kern -\theight \vbox to
    \theight{\rightline{\rlap{\box0}}%
    \vss}%
    }}%

\ShowLabelsfalse

\renewcommand{\sec}[2]{\section{#2}\label{S:#1}%
    \ifShowLabels \TeXref{{S:#1}} \fi}
\newcommand{\ssec}[2]{\subsection{#2}\label{SS:#1}%
    \ifShowLabels \TeXref{{SS:#1}} \fi}

\newcommand{\refs}[1]{Section ~\ref{S:#1}}
\newcommand{\refss}[1]{Section ~\ref{SS:#1}}

\newcommand{\reft}[1]{Theorem ~\ref{T:#1}}
\newcommand{\refl}[1]{Lemma ~\ref{L:#1}}
\newcommand{\refp}[1]{Proposition ~\ref{P:#1}}
\newcommand{\refc}[1]{Corollary ~\ref{C:#1}}
\newcommand{\refd}[1]{Definition ~\ref{D:#1}}

\newcommand{\refe}[1]{\eqref{E:#1}}

\newenvironment{thm}[1]%
    { \begin{Thm} \label{T:#1}  \ifShowLabels \TeXref{T:#1} \fi }%
    { \end{Thm} }

\renewcommand{\th}[1]{\begin{thm}{#1} \sl }
\renewcommand{\eth}{\end{thm} }

\newenvironment{lemma}[1]%
    { \begin{Lem} \label{L:#1}  \ifShowLabels \TeXref{L:#1} \fi }%
    { \end{Lem} }
\newcommand{\lem}[1]{\begin{lemma}{#1} \sl}
\newcommand{\elem}{\end{lemma}}

\newenvironment{propos}[1]%
    { \begin{Prop} \label{P:#1}  \ifShowLabels \TeXref{P:#1} \fi }%
    { \end{Prop} }
\newcommand{\prop}[1]{\begin{propos}{#1}\sl }
\newcommand{\eprop}{\end{propos}}

\newenvironment{corol}[1]%
    { \begin{Cor} \label{C:#1}  \ifShowLabels \TeXref{C:#1} \fi }%
    { \end{Cor} }
\newcommand{\cor}[1]{\begin{corol}{#1} \sl }
\newcommand{\ecor}{\end{corol}}

\newenvironment{conjec}[1]%
    { \begin{Conj} \label{Co:#1}  \ifShowLabels \TeXref{Co:#1} \fi }%
    { \end{Conj} }
\newcommand{\conj}[1]{\begin{conjec}{#1} \sl }
\newcommand{\econj}{\end{conjec}}

\newenvironment{defeni}[1]%
    { \begin{Def} \label{D:#1}  \ifShowLabels \TeXref{D:#1} \fi }%
    { \end{Def} }
\newcommand{\defe}[1]{\begin{defeni}{#1} \sl }
\newcommand{\edefe}{\end{defeni}}

\newenvironment{remark}[1]%
    { \begin{Rem} \label{R:#1}  \ifShowLabels \TeXref{R:#1} \fi }%
    { \end{Rem} }
\newcommand{\rem}[1]{\begin{remark}{#1}}
\newcommand{\erem}{\end{remark}}

\newcommand{\eq}[1]%
    { \ifShowLabels \TeXref{E:#1} \fi
       \begin{equation} \label{E:#1} }
\newcommand{\eeq}{ \end{equation} }

\newcommand{\prf}{ \begin{proof} }
\newcommand{\epr}{ \end{proof} }

\newcommand\alp{\alpha}     

\newcommand\gam{\gamma}     
     \newcommand\Del{\Delta}
\newcommand\eps{\varepsilon}

\newcommand\kap{\kappa}
\newcommand\lam{\lambda}        \newcommand\Lam{\Lambda}
\newcommand\sig{\sigma}     \newcommand\Sig{\Sigma}

\newcommand\calA{{\mathcal{A}}}

\newcommand\calD{{\mathcal{D}}}

\newcommand\calF{{\mathcal{F}}}

\newcommand\calK{{\mathcal{K}}}
\newcommand\calL{{\mathcal{L}}}
\newcommand\calM{{\mathcal{M}}}

\newcommand\calO{{\mathcal{O}}}

\newcommand\calU{{\mathcal{U}}}

\newcommand\calZ{{\mathcal{Z}}}

\newcommand\bfb{{\mathbf b}}        
        \newcommand\bfC{{\mathbf C}}
        \newcommand\bfD{{\mathbf D}}

        \newcommand\bfS{{\mathbf S}}
        
\newcommand\bfu{{\mathbf u}}        \newcommand\bfU{{\mathbf U}}

        \newcommand\bfX{{\mathbf X}}

\newcommand\RR{\mathbb{R}}

\newcommand\PP{\mathbb{P}}
\renewcommand\AA{\mathbb{A}}
\renewcommand\SS{\mathbb{S}}
\newcommand\DD{\mathbb{D}}

\newcommand\ZZ{\mathbb{Z}}

\newcommand\CC{\mathbb{C}}

 \newcommand\grg{{\mathfrak{g}}}
 \newcommand\grh{{\mathfrak{h}}}

 \newcommand\grk{{\mathfrak{k}}}
 \newcommand\grl{{\mathfrak{l}}}
 \newcommand\grm{{\mathfrak{m}}}
 \newcommand\grn{{\mathfrak{n}}}

 \newcommand\grq{{\mathfrak{q}}}

\newcommand\sdp{\times \hskip -0.3em {\raise 0.3ex
\hbox{$\scriptscriptstyle |$}}} 

\newcommand\Conv{\operatorname{Conv}}

\newcommand\End{\operatorname{End\,}}
\newcommand\Ext{\operatorname{Ext}}

\newcommand\GL{\operatorname{GL}}
\newcommand\Gr{\operatorname{Gr}}

\newcommand\Hom{\operatorname {Hom}}

\newcommand\id{\operatorname{id}}
\newcommand\Id{\operatorname{Id}}
\newcommand\im{\operatorname {im}}

\newcommand\Ker{\operatorname{Ker}}

\newcommand\Spec{\operatorname{Spec}}

\newcommand\Sym{\operatorname{Sym}}

\newcommand\uj{{\underline{j}}}

\newcommand\hatW{{\widehat{W}}}

\newcommand\x{\times}
\newcommand\ten{\otimes}

\newcommand{\ra}{\rangle}
\newcommand{\la}{\langle}

\newcommand{\nc}{\newcommand}
\nc{\renc}{\renewcommand}
\nc{\on}{\operatorname}

\nc\ol{\overline} \nc\wt{\widetilde} \nc\tboxtimes{\wt{\boxtimes}}

\renc{\SS}{{\mathbb S}} \renc{\DD}{{\mathbbD}}
\renewcommand{\AA}{{\mathbb A}}
\nc{\Fq}{{\mathbb F}_q} \nc{\Fqb}{\ol{{\mathbb F}_q}}
\nc{\Ql}{\ol{{\mathbb Q}_\ell}} \renc{\id}{\text{id}}
\nc\X{\mathcal X}

\renc{\Hom}{\on{Hom}} \nc{\Lie}{\on{Lie}} \nc{\Loc}{\on{Loc}}
\nc{\Pic}{\on{Pic}} \nc{\Bun}{\on{Bun}}
\nc{\Sh}{\on{Sh}}
\nc{\pos}{{\on{pos}}} \renc{\Conv}{\on{Conv}} \nc{\Sph}{\on{Sph}}
\renc{\Sym}{\on{Sym}}
\nc{\BunBb}{\overline{\Bun}_B} \nc{\Buno}{\overset{o}{\Bun}}
\nc{\BunPb}{{\overline{\Bun}_P}}
\nc{\BunBM}{\overline{\Bun}_{B(M)}}
\nc{\BunPbw}{{\widetilde{\Bun}_P}}
\nc{\BunBP}{\widetilde{\Bun}_{B,P}} \nc{\GUb}{\overline{G/U}}
\nc{\GUPb}{\overline{G/U(P)}}

\nc{\Hhom}{\underline{\on{Hom}}} \nc\syminfty{\on{Sym}^{\infty}}
\nc\lal{\ol{\lambda}} \nc\xl{\ol{x}} \nc\thl{\ol{\theta}}
\nc\nul{\ol{\nu}} \nc\mul{\ol{\mu}} \nc\Sum\Sigma
\nc{\hl}{\overset{\leftarrow}h}
\nc{\hr}{\overset{\rightarrow}h} \nc{\M}{{\mathcal M}}
\nc{\N}{{\mathcal N}} \nc{\F}{{\mathcal F}} \nc{\D}{{\mathcal D}}
\nc{\Q}{{\mathcal Q}} \nc{\Y}{{\mathcal Y}} \nc{\G}{{\mathcal G}}
\nc{\E}{{\mathcal E}} \nc{\CalC}{{\mathcal C}}
\nc\Dh{\widehat{\D}}

\nc{\C}{{\mathcal C}} \nc{\K}{{\mathcal K}}
\renewcommand{\H}{{\mathcal H}}

\nc{\T}{{\mathcal T}} \nc{\V}{{\mathcal V}} \renc{\P}{{\mathcal
P}} \nc{\A}{{\mathcal A}} \nc{\B}{{\mathcal B}} \nc{\U}{{\mathcal
U}}

\renc{\Gr}{\on{Gr}}

\nc{\frn}{{\check{\mathfrak u}(P)}} \nc{\p}{\mathfrak p}
\nc{\q}{\mathfrak q} \nc\f{{\mathfrak f}}

\nc{\qo}{{\mathfrak q}} \nc{\po}{{\mathfrak p}} \nc{\s}{{\mathfrak
s}} \nc\w{\text{w}}

\nc\mathi\iota \renc\Spec{\on{Spec}} \nc\Mod{\on{Mod}}
\nc{\tw}{\widetilde{\mathfrak t}} \nc{\pw}{\widetilde{\mathfrak
p}} \nc{\qw}{\widetilde{\mathfrak q}} \nc{\jw}{\widetilde j}

\nc{\I}{\mathcal I}

\nc{\lambdach}{{\check\lambda}} \nc{\Lambdach}{{\check\Lambda}{}}
\nc{\much}{{\check\mu}} \nc{\omegach}{{\check\omega}}
\nc{\nuch}{{\check\nu}} \nc{\etach}{{\check\eta}}
\nc{\alphach}{{\check\alpha}} \nc{\betach}{{\check\beta}}
\nc{\rhoch}{{\check\rho}} \nc{\ch}{{\check h}}

\nc{\Hb}{\overline{\H}}

\nc{\BA}{{\mathbb{A}}} \nc{\BC}{{\mathbb{C}}}
\nc{\BG}{{\mathbb{G}}} \nc{\BM}{{\mathbb{M}}}
\nc{\BN}{{\mathbb{N}}} \nc{\BP}{{\mathbb{P}}}
\nc{\BR}{{\mathbb{R}}} \nc{\BZ}{{\mathbb{Z}}}
\nc{\BS}{{\mathbb{S}}}

\nc{\CA}{{\mathcal{A}}} \nc{\CB}{{\mathcal{B}}}

\nc{\CE}{{\mathcal{E}}} \nc{\CF}{{\mathcal{F}}}
\nc{\CG}{{\mathcal{G}}} \nc{\CL}{{\mathcal{L}}}
\nc{\CM}{{\mathcal{M}}} \nc{\CN}{{\mathcal{N}}}
\nc{\CK}{{\mathcal{K}}} \nc{\CO}{{\mathcal{O}}}
\nc{\CP}{{\mathcal{P}}} \nc{\CQ}{{\mathcal{Q}}}
\nc{\CR}{{\mathcal{R}}} \nc{\CS}{{\mathcal{S}}}
\nc{\CT}{{\mathcal{T}}} \nc{\CU}{{\mathcal{U}}}
\nc{\CV}{{\mathcal{V}}} \nc{\CW}{{\mathcal{W}}}
\nc{\CZ}{{\mathcal{Z}}}

\nc{\cM}{{\check{\mathcal M}}{}} \nc{\csM}{{\check{\mathcal A}}{}}
\nc{\obM}{{\overset{\circ}{\mathbf M}}{}}
\nc{\oCA}{{\overset{\circ}{\mathcal A}}{}}
\nc{\obA}{{\overset{\circ}{\mathbf A}}{}}
\nc{\ooM}{{\overset{\circ}{M}}{}}
\nc{\osM}{{\overset{\circ}{\mathsf M}}{}}
\nc{\vM}{{\overset{\bullet}{\mathcal M}}{}}
\nc{\nM}{{\underset{\bullet}{\mathcal M}}{}}
\nc{\obD}{{\overset{\circ}{\mathbf D}}{}}
\nc{\cp}{{\overset{\circ}{\mathbf p}}{}}
\nc{\ofZ}{{\overset{\circ}{\mathfrak Z}}{}}

\nc{\fa}{{\mathfrak{a}}} \nc{\fb}{{\mathfrak{b}}}
\nc{\fg}{{\mathfrak{g}}} \nc{\fgl}{{\mathfrak{gl}}}
\nc{\fh}{{\mathfrak{h}}} \nc{\fj}{{\mathfrak{j}}}
\nc{\fm}{{\mathfrak{m}}} \nc{\fn}{{\mathfrak{n}}}
\nc{\fu}{{\mathfrak{u}}} \nc{\fp}{{\mathfrak{p}}}
\nc{\fr}{{\mathfrak{r}}} \nc{\fs}{{\mathfrak{s}}}
\nc{\fsl}{{\mathfrak{sl}}} \nc{\hsl}{{\widehat{\mathfrak{sl}}}}
\nc{\hgl}{{\widehat{\mathfrak{gl}}}}
\nc{\hg}{{\widehat{\mathfrak{g}}}}
\nc{\chg}{{\widehat{\mathfrak{g}}}{}^\vee}
\nc{\hn}{{\widehat{\mathfrak{n}}}}
\nc{\chn}{{\widehat{\mathfrak{n}}}{}^\vee}

\nc{\fA}{{\mathfrak{A}}} \nc{\fB}{{\mathfrak{B}}}
\nc{\fD}{{\mathfrak{D}}} \nc{\fE}{{\mathfrak{E}}}
\nc{\fF}{{\mathfrak{F}}} \nc{\fG}{{\mathfrak{G}}}
\nc{\fK}{{\mathfrak{K}}} \nc{\fL}{{\mathfrak{L}}}
\nc{\fM}{{\mathfrak{M}}} \nc{\fN}{{\mathfrak{N}}}
\nc{\fP}{{\mathfrak{P}}} \nc{\fU}{{\mathfrak{U}}}
\nc{\fV}{{\mathfrak{V}}} \nc{\fZ}{{\mathfrak{Z}}}

\nc{\bb}{{\mathbf{b}}} \nc{\bc}{{\mathbf{c}}}
\nc{\bd}{{\mathbf{d}}} \nc{\be}{{\mathbf{e}}}
\nc{\bj}{{\mathbf{j}}} \nc{\bn}{{\mathbf{n}}}
\nc{\bp}{{\mathbf{p}}} \nc{\bq}{{\mathbf{q}}}
\nc{\bu}{{\mathbf{u}}} \nc{\bv}{{\mathbf{v}}}
\nc{\bx}{{\mathbf{x}}} \nc{\bs}{{\mathbf{s}}}
\nc{\by}{{\mathbf{y}}} \nc{\bw}{{\mathbf{w}}}
\nc{\bA}{{\mathbf{A}}} \nc{\bK}{{\mathbf{K}}}
\nc{\bB}{{\mathbf{B}}} \nc{\bC}{{\mathbf{C}}}
\nc{\bD}{{\mathbf{D}}} \nc{\bH}{{\mathbf{H}}}
\nc{\bM}{{\mathbf{M}}} \nc{\bN}{{\mathbf{N}}}
\nc{\bV}{{\mathbf{V}}} \nc{\bW}{{\mathbf{W}}}
\nc{\bX}{{\mathbf{X}}} \nc{\bZ}{{\mathbf{Z}}}
\nc{\bS}{{\mathbf{S}}}

\nc{\sA}{{\mathsf{A}}} \nc{\sB}{{\mathsf{B}}}
\nc{\sC}{{\mathsf{C}}} \nc{\sD}{{\mathsf{D}}}
\nc{\sF}{{\mathsf{F}}} \nc{\sK}{{\mathsf{K}}}
\nc{\sM}{{\mathsf{M}}} \nc{\sO}{{\mathsf{O}}}
\nc{\sQ}{{\mathsf{Q}}} \nc{\sP}{{\mathsf{P}}}
\nc{\sZ}{{\mathsf{Z}}} \nc{\sfp}{{\mathsf{p}}}
\nc{\sr}{{\mathsf{r}}} \nc{\sg}{{\mathsf{g}}}
\nc{\sff}{{\mathsf{f}}} \nc{\sfb}{{\mathsf{b}}}
\nc{\sfc}{{\mathsf{c}}} \nc{\sd}{{\mathsf{d}}}

\nc{\BK}{{\bar{K}}}

\nc{\tA}{{\widetilde{\mathbf{A}}}}
\nc{\tB}{{\widetilde{\mathcal{B}}}}
\nc{\tg}{{\widetilde{\mathfrak{g}}}} \nc{\tG}{{\widetilde{G}}}
\nc{\TM}{{\widetilde{\mathbb{M}}}{}}
\nc{\tO}{{\widetilde{\mathsf{O}}}{}}
\nc{\tU}{{\widetilde{\mathfrak{U}}}{}} \nc{\TZ}{{\tilde{Z}}}
\nc{\tx}{{\tilde{x}}} \nc{\tbv}{{\tilde{\bv}}}
\nc{\tfP}{{\widetilde{\mathfrak{P}}}{}} \nc{\tz}{{\tilde{\zeta}}}
\nc{\tmu}{{\tilde{\mu}}}

\nc{\urho}{\underline{\rho}} \nc{\uB}{\underline{B}}
\nc{\uC}{{\underline{\mathbb{C}}}} \nc{\ui}{\underline{i}}
\renc{\uj}{\underline{j}} \nc{\ofP}{{\overline{\mathfrak{P}}}}

\renc{\eps}{\varepsilon} \nc{\hrho}{{\hat{\rho}}}

\nc{\one}{{\mathbf{1}}} \nc{\two}{{\mathbf{t}}}

\nc{\Rep}{{\mathop{\operatorname{\rm Rep}}}}
\nc{\Tot}{{\mathop{\operatorname{\rm Tot}}}}
\renc{\Ker}{{\mathop{\operatorname{\rm Ker}}}}
\nc{\Hilb}{{\mathop{\operatorname{\rm Hilb}}}}
\renc{\End}{{\mathop{\operatorname{\rm End}}}}
\renc{\Ext}{{\mathop{\operatorname{\rm Ext}}}}
\nc{\CHom}{{\mathop{\operatorname{{\mathcal{H}}\it om}}}}
\renc{\GL}{{\mathop{\operatorname{\rm GL}}}}
\nc{\gr}{{\mathop{\operatorname{\rm gr}}}}
\renc{\Id}{{\mathop{\operatorname{\rm Id}}}}
\nc{\de}{{\mathop{\operatorname{\rm def}}}}
\nc{\length}{{\mathop{\operatorname{\rm length}}}}

\nc{\Cliff}{{\mathsf{Cliff}}}
\nc{\Fl}{\on{Fl}} \nc{\Fib}{{\mathsf{Fib}}}
\nc{\Coh}{{\mathsf{Coh}}} \nc{\FCoh}{{\mathsf{FCoh}}}

\nc{\reg}{{\text{\rm reg}}}

\nc{\cplus}{{\mathbf{C}_+}} \nc{\cminus}{{\mathbf{C}_-}}
\nc{\cthree}{{\mathbf{C}_*}} \nc{\Qbar}{{\bar{Q}}}

\nc{\bh}{{\bar{h}}} \nc{\bOmega}{{\overline{\Omega}}}

\nc{\seq}[1]{\stackrel{#1}{\sim}} \nc\QM{{\mathcal {QM}}}

\nc{\chH}{\check H} \nc{\chM}{\check M} \nc{\aff}{{\on{aff}}}

\nc{\chh}{\check \grh}

\renewcommand\chn{\check \grn}

\renewcommand\chg{\check{\grg}}

\nc\chT{\check T}
\renewcommand\Ext{\operatorname{Ext}}
\renewcommand\d{\partial}
\newcommand\chlam{\check \lam}

\begin{document}
Dedicated to A.~Joseph on the occasion of his 60th birthday

\bigskip

\title[Instanton counting via affine Lie algebras II]{Instanton counting via affine Lie algebras II:
from Whittaker vectors to the Seiberg-Witten prepotential}
\author{Alexander Braverman and Pavel Etingof}
\email{braval@math.brown.edu; etingof@math.mit.edu}
\address{Department of Mathematics, Brown
University, 151
Thayer street, Providence RI, USA and Einstein Institute of Mathematics,
Edmond J. Safra Campus, Givat Ram
The Hebrew University of Jerusalem
Jerusalem, 91904, Israel}
\address{Department of Mathematics, Massachusetts Institute of
  Technology, 77 Mass. Ave., Cambridge MA, 02139, USA}

\maketitle

\begin{abstract}
Let $G$ be a simple
simply connected algebraic group over $\CC$
with Lie algebra $\grg$. Given a parabolic subgroup $P\subset G$, in
\cite{Br} the first author introduced a certain generating function
$Z_{G,P}^{\aff}$. Roughly speaking, these functions count (in a
certain sense) framed $G$-bundles on $\PP^2$ together with a
$P$-structure on a fixed (horizontal) line in $\PP^2$. When $P=B$ is
a Borel subgroup, the function $Z_{G,B}^{\aff}$ was identified in
\cite{Br} with the Whittaker matrix coefficient in the universal
Verma module over the  affine Lie algebra $\chg_{\aff}$ (here we
denote by $\grg_{\aff}$ the affinization of $\grg$ and by
$\chg_{\aff}$ the Lie algebra whose root system is dual to that of
$\grg_{\aff}$).

For $P=G$ (in this case we shall write $\calZ_G^{\aff}$ instead of
$\calZ_{G,P}^{\aff}$) and $G=SL(n)$ the above generating function
was introduced by Nekrasov (cf. \cite{nek}) and studied thoroughly
in \cite{nayo} and \cite{neok}. In particular, it is shown in {\it
loc. cit.} that the leading term of certain asymptotic of
$\calZ_G^{\aff}$ is given by the (instanton part of the) {\it
Seiberg-Witten prepotential} (for $G=SL(n)$). The prepotential is
defined using the geometry of the (classical) periodic Toda
integrable system. This result was conjectured in \cite{nek}.

The purpose of this paper is to extend these results to arbitrary
$G$. Namely, we use the above description of the function
$\calZ_{G,B}^{\aff}$ to show that the leading term of its
asymptotic (similar to the one studied in \cite{nek} for $P=G$) is
given by the instanton part of the prepotential constructed via the
Toda system attached to the Lie algebra $\chg_{\aff}$. This part is
completely algebraic and does not use the original algebro-geometric
definition of $\calZ_{G,B}^{\aff}$. We then show that for fixed $G$
these asymptotic are the same for {\it all} functions
$\calZ_{G,P}^{\aff}$.
\end{abstract}

\sec{}{Introduction} \ssec{}{The partition function} This paper
has grown out of a (still unsuccessful) attempt to understand the
following object. Let $K$ be a simple
\footnote{In this paper by a simple Lie (or algebraic) group we mean a group
whose Lie algebra is simple}
simply connected compact Lie
group and let $d$ be a non-negative integer. Denote by $\calM_K^d$
the moduli space of (framed) $K$-instantons on $\RR^4$  of second
Chern class $-d$. This space can be naturally embedded into a
larger {\it Uhlenbeck space} $\calU_K^d$. Both spaces admit a
natural action of the group $K$ (by changing the framing at
$\infty$) and the torus $(S^1)^2$ acting on $\RR^4$ after choosing
an identification $\RR^4\simeq \CC^2$. Moreover,  the maximal
torus of $K\x (S^1)^2$ has unique fixed point on $\calU_K^d$. Thus
we may consider (cf. \cite{Br}, \cite{nayo} or \cite{nek}  for
precise definitions) the {\it equivariant integral}
$$
\int\limits_{\calU^d_K}1^d
$$
of the unit $K\times (S^1)^2$-equivariant cohomology class (which we
denote by $1^d$) over
$\calU_K^d$; the integral takes values in the field $\calK$
of fractions of the algebra $\calA=H^*_{K\x
(S^1)^2}(pt)$ \footnote{In this paper we always consider
cohomology with complex coefficients.}. Note that $\calA$ is
canonically isomorphic to the algebra of polynomial functions on
$\grk\x\RR^2$ (here $\grk$ denotes the Lie algebra of $K$) which
are invariant with respect to the adjoint action of $K$ on $\grk$.
Thus each $\int\limits_{\calU^d_K}1^d$ may be naturally regarded as a
rational function of $a\in\grk$ and $(\eps_1,\eps_2)\in \RR^2$.

Consider now the generating function
$$
\calZ=\sum\limits_{d=0}^\infty Q^d \int\limits_{\calU_K^d} 1^d.
$$
It can (and should) be thought of as a function of the variables
$Q$ and $a,\eps_1,\eps_2$ as before. In \cite{nek} it was
conjectured that the first term of the asymptotic in the limit
$\underset{\eps_1,\eps_2\to 0}\lim \ln \calZ$ is closely related to
the {\it Seiberg-Witten prepotential} of $K$. For $K=SU(n)$ this
conjecture has been proved in \cite{neok} and \cite{nayo}. Also in
\cite{nek} an explicit combinatorial expression for $\calZ$ has
been found.

\ssec{}{Algebraic version}In \cite{Br} the first author has defined some more
general partition functions containing the function $\calZ_K$ as a special case.
Let us recall that definition. First of all, it will be
convenient for us to make the whole situation completely
algebraic.

Namely, let $G$ be a complex simple algebraic group whose
maximal compact subgroup is isomorphic to $K$. We shall denote by
$\grg$ its Lie algebra. Let also $\bfS=\PP^2$ and denote by
$\bfD_\infty\subset \bfS$ the "straight line at $\infty$"; thus
$\bfS\backslash\bfD_\infty=\AA^2$. It is well-known that $\calM^d_K$
is isomorphic to the moduli space $\Bun_G^d(\bfS,\bfD_\infty)$ of
principal $G$-bundles on $\bfS$ of second Chern class $-d$
endowed with a trivialization on
$\bfD_\infty$. When it does not lead to
a confusion we shall write $\Bun_G$ instead of
$\Bun_G(\bfS,\bfD_\infty)$. The algebraic analog of $\calU^d_K$ has
been constructed in \cite{bfg}; we denote this algebraic variety
by $\calU_G^d$. This variety is endowed with a natural action on
$G\x (\CC^*)^2$.
\ssec{}{Parabolic generalization of the partition function}Let
$\bfC\subset\bfS$ denote the standard horizontal line. Choose a parabolic
subgroup $P\subset G$. Let $\Bun_{G,P}$ denote the moduli space of
the following objects:

1) A principal $G$-bundle $\calF_G$ on $\bfS$;

2) A trivialization of $\calF_G$ on $\bfD_\infty\subset \bfS$;

3) A reduction of $\calF_G$ to $P$ on $\bfC$ compatible with the
trivialization of $\calF_G$ on $\bfC\cap \bfD_{\infty}$.

Let us describe the connected components of $\Bun_{G,P}$. Let $M$
be the Levi group of $P$. Denote by $\chM$ the {\it Langlands
dual} group of $M$ and let $Z(\chM)$ be its center. We denote by
$\Lam_{G,P}$ the lattice of characters of $Z(\chM)$. Let also
$\Lam_{G,P}^{\aff}=\Lam_{G,P}\x \ZZ$ be the lattice of characters
of $Z(\chM)\x\CC^*$. Note that $\Lam_{G,G}^{\aff}=\ZZ$.

The lattice $\Lam^{\aff}_{G,P}$ contains a canonical semi-group
$\Lam^{\aff,\pos}_{G,P}$ of positive elements (cf. \cite{bfg} and \cite{Br}).
It is not difficult to see that the connected
components of $\Bun_{G,P}$ are parameterized  by the elements of
$\Lam_{G,P}^{\aff,\pos}$:
$$
    \Bun_{G,P}=\bigcup\limits_{\theta_{\aff}\in\Lam_{G,P}^{\aff,\pos}}
    \Bun_{G,P}^{\theta_{\aff}}.
$$

Typically, for $\theta_{\aff}\in \Lam_{G,P}^{\aff}$ we shall write
$\theta_\aff=(d,\theta)$ where $\theta\in \Lam_{G,P}$ and $d\in \ZZ$.

Each $\Bun_{G,P}^{\theta_{\aff}}$ is naturally acted on by $P\x(\CC^*)^2$;
by embedding $M$ into $P$ we get an action of $M\x(\CC^*)^2$ on
$\Bun_{G,P}^{\theta_{\aff}}$. In \cite{bfg} we define for each
$\theta_{\aff}\in\Lam_{G,P}^{\aff,\pos}$ a certain Uhlenbeck scheme
$\calU_{G,P}^{\theta_{\aff}}$ which contains $\Bun_{G,P}^{\theta_{\aff}}$ as a dense
open subset. The scheme $\calU_{G,P}^{\theta_{\aff}}$ still admits an
action of $M\x(\CC^*)^2$.

We want to do some equivariant intersection theory on the spaces
$\calU_{G,P}^{\theta_{\aff}}$. For this let us denote by $\calA_{M\x(\CC^*)^2}$
the algebra $H^*_{M\x (\CC^*)^2}(pt,\CC)$. Of course this is just
the algebra of $M$-invariant polynomials on $\grm\x\CC^2$.
Let also $\calK_{M\x(\CC^*)^2}$ be its field of fractions. We can think
about elements of $\calK_{M\x (\CC^*)^2}$ as rational functions on
$\grm\x\CC^2$ which are invariant with respect to the adjoint action.

Let $T\subset M$ be a maximal torus. Then one can show that
$(\calU_{G,P}^{\theta_{\aff}})^{T\x(\CC^*)^2}$ consists of one point.
This guarantees that we may consider the integral
$\int\limits_{\calU_{G,P}^{\theta_{\aff}}}1_{G,P}^{\theta_{\aff}}$ where
$1_{G,P}^{\theta_{\aff}}$
denotes the unit class in
$H^*_{M\x(\CC^*)^2}(\calU_{G,P}^{\theta_{\aff}},\CC)$.
The result can be thought
of as a rational function on $\grm\x\CC^2$ which is invariant with
respect to the adjoint action of $M$. Define
\eq{partition-affine}
    \calZ_{G,P}^{\aff}=\sum\limits_{\theta\in\Lam_{G,P}^{\aff}} \grq_{\aff}^{\theta_{\aff}}
    \ \int\limits_{\calU_{G,P}^{\theta_{\aff}}}1_{G,P}^{\theta_{\aff}}
\end{equation}
(we refer the reader to Section 2 of \cite{Br} for a detailed discussion of equivariant integration).
One should think of $\calZ_{G,P}^{\aff}$ as a formal power series
in $\grq_{\aff}\in Z(\chM)\x\CC^*$ with values in the space of ad-invariant
rational functions on $\grm\x\CC^2$. Typically, we shall write
$\grq_{\aff}=(\grq,Q)$ where $\grq\in Z(\chM)$ and $Q\in\CC^*$.
Also we shall denote an element of $\grm\x\CC^2$ by
$(a,\eps_1,\eps_2)$ (note that for general $P$ (unlike in the case $P=G$)
the function $\calZ_{G,P}^{\aff}$ is not symmetric with respect to
switching
$\eps_1$ and $\eps_2$).
Here is the main result of this paper:

\th{main}
Let $P\subset G$ be a parabolic subgroup as a above.
\begin{enumerate}
\item
There exists
a function $\calF^{inst}\in \CC(a)[[Q]]$ such that
\eq{limit}
\underset{\eps_1\to 0}\lim\underset{\eps_2\to 0}\lim \, \eps_1\eps_2\ln Z_{G,P}^{\aff}=\calF^{inst}(a,Q).
\end{equation}
In particular, the above limit does not depend on $\grq$ and it is the same for all
$P$.
\item
The function $\calF^{inst}(a,Q)$ is equal to {\it the instanton part of the
Seiberg-Witten prepotential} of the affine Toda system associated with the
Langlands dual Lie algebra $\chg_{\aff}$ (cf. \refs{highdim} for the explanation of these words).
\end{enumerate}
\eth Since the function $\calZ_G^{\aff}$ is symmetric in $\eps_1$
and $\eps_2$, \reft{main} implies the following result:
\cor{nikita} The function $\eps_1\eps_2\ln \calZ_G^{\aff}$ is
regular when both $\eps_1$ and $\eps_2$ are set to $0$. Moreover,
one has
$$
(\eps_1\eps_2\ln \calZ_G^{\aff})|_{\eps_1=\eps_2=0}=\calF^{inst}.
$$

\ecor
\refc{nikita} was conjectured by N.~Nekrasov in \cite{nek} (in fact \cite{nek}
contains only the formulation for $G=SL(n)$ but the generalization to other groups
is straightforward). For $G=SL(n)$ Nekrasov's conjecture was proved in \cite{nayo} and \cite{neok}.
Also, more recently, this conjecture was proved in \cite{nesh} for all classical groups.
These papers, however, utilize methods which are totally different from ours. In particular, in our
approach the existence of the partition functions $\calZ_{G,P}^{\aff}$ for $P\neq G$
(in particular, for $P$ being the Borel subgroup) plays a crucial
role.

In fact, we are going to prove the following slightly stronger version
of \reft{main}:
\th{main'}
\begin{enumerate}
\item
\reft{main} holds for $P=B$.
\item
For every parabolic subgroup $P\subset G$ one has
$$
\underset{\eps_2\to 0}\lim \eps_2(\ln \calZ_{G,P}^{\aff} - \ln
\calZ_{G}^{\aff})=0.
$$
\end{enumerate}
\eth

\ssec{}{Plan of the proof}
Let us explain the idea of the proof of \reft{main'}.
The second part is in fact rather routine so let us explain the idea of the proof of the first
part.

The "Borel" partition function $\calZ_{G,B}^{\aff}$  was realized
in \cite{Br} as the Whittaker matrix coefficient in the universal
Verma module over the Lie algebra $\chg_{\aff}$. As a corollary
one gets that the function $\calZ_{G,B}^{\aff}$ is an
eigenfunction of {\it the non-stationary Toda hamiltonian}
associated with the affine Lie algebra $\chg_{\aff}$ (cf.
Corollary 3.7 from \cite{Br} for the precise statement; we use
\cite{etingof} as our main reference about Toda hamiltonians).

It turns out that this is all that we have to use in order to prove \reft{main'}(1). Namely, in this paper
(cf. \refs{1d} and \refs{highdim}) we introduce the notion of the
Seiberg-Witten prepotential (more precisely, its instanton part)
for a very general class of non-stationary Schr\"odinger operators in such a way that by the definition it is equal to
some asymptotic of the (in some sense) universal eigenfunction of this operator
(we were unable to find such a definition in the
literature).
Usually the prepotential is attached to a classical completely integrable system
(our main references on the definition of the Seiberg-Witten prepotential are \cite{NaYo1} and
\cite{nek}). We show that in the integrable case our definition of the prepotential coincides
with the one from {\it loc. cit}.

\ssec{}{Acknowledgments} A.B. would like to thank
M.~Finkelberg, A.~Gorsky, A.~Marshakov,
A.~Mironov, I.~Krichever, H.~Nakajima, N.~Nekrasov, A.~Okounkov and
Y.~Oz for very valuable discussions on the subject. The work of A.B.
was partially supported by the NSF grant DMS-0300271. The work of
P.E. was partially supported by the NSF grant DMS-9988796 and the
CRDF grant RM1-2545-MO-03.

\sec{1d}{Schr\"odinger operators and the prepotential: the
one-dimensional case}

\ssec{}{Schr\"odinger operators}Let $x\in \CC$ and let $U(x)$ be a
trigonometric polynomial in $x$ - i.e. a polynomial in $e^x$ and
$e^{-x}$. Let also $\hbar$ and $Q$ be formal variables. We want to
study the eigenvalues of the Schr\"odinger operator
$$
T=\hbar^2\frac{d^2}{dx^2}+QU(x).
$$
More precisely,  for each $a\in  \CC$ let
$W_a$ denote the space $e^{\frac{ax}{\hbar}}\CC(\hbar)[e^x,e^{-x}][[Q]]$ with the natural
action of the algebra of linear differential operators in $x$.
Then we would like to look for eigenfunctions of $T$ in $W_a$.
After conjugating $T$ with $e^{\frac{ax}{\hbar}}$ the operator $T$
turns into the operator
$$
\hbar^2\frac{d^2}{dx^2}+2\hbar a\frac{d}{dx}+QU(x)+a^2.
$$
Let
$$
T^a=\hbar^2\frac{d^2}{dx^2}+2\hbar a\frac{d}{dx}+QU(x).
$$
We now want to look for eigenfunctions of $T^a$ in $W_0$
(this problem is obviously equivalent to finding eigenfunctions
of $T$ in $W_a$). In fact,
we want them to depend nicely on $a$, so we set
$W=\CC(a,\hbar)[e^x,e^{-x}][[Q]]$ and we want to look for eigenfunctions of
$T^a$ (considered now as a differential operator with coefficients in
$\CC(a,\hbar)[[Q]]$) in $W$.
\prop{stat}
\begin{enumerate}
\item
There exist $\psi\in W$ and $b\in Q\CC(a,\hbar)[[Q]]$ such that
\eq{eig}
T^a\psi=b\psi
\end{equation}
and such that $\psi=1+O(Q)$.
Moreover, under such conditions $b$ is unique  and
$\psi$ is unique up to multiplication by an element of
$1+Q\CC(a,\hbar)[[Q]]$.
\item
Let $\phi=\hbar\ln \psi$ (note that $\phi$ is defined uniquely up
to adding an element of $Q\CC(a,\hbar)[[Q]])$. Then $\phi$ is regular
at $\hbar=0$ provided this is true for its constant term.
\footnote{By the "constant term" we shall always mean the constant term of a trigonometric polynomial.
The reader should not confuse this with the notion of
"free term" by which we always mean the coefficient
of the $0$-th power of the variable in a formal power series.}
\item
The limit $v(a,Q):=\underset{\hbar\to 0}\lim b(a,\hbar,Q)$ exists in
$\CC(a)[[Q]]$.
\end{enumerate}
\eprop
\prf

Let us prove the first assertion.
Let us write
$$
\psi=\sum\limits_{n=0}^{\infty}\psi_nQ^n\qquad{\text{and}}\qquad
b=\sum\limits_{n=0}^{\infty}b_n Q^n.
$$
Note that $\psi_0=1$ and thus automatically $b_0=0$. Thus the
equation
\refe{eig} becomes
\eq{rew}
\hbar^2 \psi_n''+2\hbar a
\psi_n'+U(x)\psi_{n-1}=\sum\limits_{i=0}^{n-1} b_{n-i}\psi_{i}.
\end{equation}
which should be valid for each $n>0$ (here and in what follows
the prime denotes the derivative of a function with respect to $x$).
It is enough for us to prove that
the system of equations \refe{rew} has a unique solution if we require that
for all $n>0$ the constant term of the function $\psi_n$ is equal to 0.

Equation \refe{rew} is equivalent to
\eq{rew'}
\hbar^2\psi_n''+2\hbar a
\psi_n'=-U(x)\psi_{n-1}+\sum\limits_{i=0}^{n-1} b_{n-i}\psi_{i}.
\end{equation}
where the left hand side is just the differential operator
$D=\hbar\frac{d^2}{dx^2}+2\hbar a\frac{d}{dx}$ applied to $\psi_n$, and
the right hand side only depends on the $\psi_i$'s with $i<n$.

Let us now argue by induction on $n$. By the induction hypothesis we
assume that $\psi_i$ and $b_i$ have already been uniquely determined
for all $i<n$. Note that the operator $D$ has the following properties
(whose verification is left to the reader):

1) $\ker D$ consists of constant (i.e. independent of $x$) functions.

2) $\im D$ consists of all trigonometric polynomials whose constant
   term is equal to $0$.

Observe now that the coefficient of $b_n$ in the RHS of \refe{rew'} is $\psi_0=1$.
Thus property 2) above determines $b_n$ uniquely -- it has to be chosen
so that the constant term of the RHS is equal to $0$. If $b_n$ is
chosen in this way then there exists some $\psi_n$ satisfying
\refe{rew'}. {\it A priori} such $\psi_n$ is defined uniquely up to
adding a constant trigonometric polynomial, but the requirement
that the constant term of $\psi_n$ is equal to $0$ determines $\psi_n$ uniquely.

Let us prove the second and third assertions (this is a standard WKB argument
which we include for the sake of completeness). Let us write
$$
\phi=\hbar\ln \psi
$$
(the logarithm is taken in the sense of formal power series in $Q$; this makes sense because
$\psi_0=1$).

Let us rewrite \refe{eig} in terms of $\phi$.
We get
\eq{phi}
(\phi')^2+\hbar\phi''+2a\phi'+Q U(x)=b.
\end{equation}
Let us now look for a solution $\phi$ of the form
$$
\phi=\sum\limits_{n=1}^{\infty}\phi_n Q^n
$$

Then \refe{phi} is equivalent to the following system of equations:
\eq{phione}
\hbar\phi_1''+2a\phi_1'=b_1-U(x)
\end{equation}
and
\eq{phin}
\hbar\phi_n''+2a\phi_n'=b_n-\sum\limits_{i=1}^{n-1} \phi_i'\phi_{n-i}'
\end{equation}
for all $n>1$.

Without loss of generality we may assume that the constant term of all $\phi_n$ is equal to zero.
We need to show that under such conditions all $\phi_n$ and $b_n$
are regular when $\hbar=0$.
Let us prove by induction in $k$ that
the statement is valid for $n\le k$.
If $k=0$, the statement is clear, so let $k>0$; we need to prove
the statement for $n=k$.
By the induction assumption
 we may assume that $\sum\limits_{i=1}^{n-1} \phi_i'\phi_{n-i}'$
is regular at $\hbar=0$. Arguing as before we see that if \refe{phin} has a solution then $b_n$ has to be equal to the
constant term of $\sum\limits_{i=1}^{n-1} \phi_i'\phi_{n-i}'$
for $n>1$ and of $U(x)$ for $n=1$, and thus it is also regular at
$\hbar=0$. Thus the right hand side
of \refe{phin} is regular at $\hbar=0$. This immediately implies
that the same is true for $\phi_n$.
\epr
\ssec{exp1d}{Explicit calculation of $\lim\limits_{\hbar\to 0} b$ via periods}
We now want to explain how to evaluate the
function $\underset{\hbar\to 0}\lim  \, b(a,\hbar,Q)=v(a,Q)$ using
period integrals on a certain algebraic curve. More precisely, we
are going to express $a$ as a (multi-valued) function of $v$ and $Q$
which will be written in terms of such periods. Let ${\varphi}$
denote the limit of $\phi$ as $\hbar\to 0$.  Then we have the equation \eq{phi0}
({\varphi}\, ')^2+2a{\varphi}\, '=v-Q U(x).
\end{equation}
In other words, ${\varphi}\, '$ satisfies a quadratic equation.
Thus we may write
$$
{\varphi}\, '=-a+\sqrt{a^2+v-QU(x)}.
$$
This is an equality of formal power series in $Q$. The square root is chosen in such a way that
the right hand side is equal to $0$ when $Q=0$ (note we automatically have $v=0$ when $Q=0$).

Recall, however that ${\varphi}$ was a trigonometric polynomial.
This implies that
$$
\int\limits_0^{2\pi i} {\varphi}\, 'dx=0.
$$
This is equivalent to the equation
\eq{brbrbr}
2\pi i a=\int\limits_0^{2\pi i} \sqrt{a^2+v-QU(x)}dx
\end{equation}

Set $w=e^x$ and recall that $U(x)=P(w)$ for some polynomial $P$ in $w$ and $w^{-1}$.
Set also $u=a^2+v$ and consider the algebraic curve $C=C_u$ which is the projectivization of
the affine curve given by the equation
$$
z^2+QP(w)=u.
$$

We claim that we may write $a$ locally as a function $a(u,Q)$ of $u$ and $Q$. Namely, first of
all $a_0:=a(u,0)$ must satisfy $a_0^2=u$. Let us locally choose
one of the square roots. Then the function $a$ is found as a
series $a_0+a_1Q+...+a_nQ^n+...$, where $a_i$ with $i>0$ are
found recursively.

Note that when $Q=0$ the above curve breaks into two components corresponding
to $z=\pm a_0$.
Let $A=A_{u,Q}$ denote the one-dimensional cycle in $C$ satisfying the following conditions:

1) The projection of $A$ to the $w$-plane is an isomorphism between $A$ and the unit circle.

2) $A$ depends continuously on $Q$ and when $Q=0$ it lies in the component of $C$ corresponding
to $z=a_0$.

Such a cycle is unique at least for small values of $Q$. Thus the
equation \refe{brbrbr} becomes equivalent to
\eq{periods}
a=\frac{1}{2\pi i}\oint\limits_A z\frac{dw}{w}.
\end{equation}
Note that $z\frac{dw}{w}$ is a well-defined meromorphic differential on $C$.
Note also that $C$ and $A$ depend only on $Q$ and $u$; thus we may think of \refe{periods} as expressing
$a$ as a function of $u$ and $Q$.
\ssec{}{Eigenfunctions of non-stationary Schr\"odinger operators}
Let us now change our problem a little. Introduce one more
variable $\kap$ and define new operators
$$
\calL=T-\kap Q\frac{\d}{\d Q}\qquad\text{and}\qquad
\calL^a=T^a-\kap Q\frac{\d}{\d Q}.
$$
Let us now look for solutions of the equation
\eq{non-st}
\calL^a\Psi=0
\end{equation}
where $\Psi\in \CC(a,\hbar,\kap)[e^x,e^{-x}][[Q]]$ (we shall denote this space by $W(\kap)$).
Of course this equation is equivalent
to the equation
$$
\calL(e^{\frac{ax}{\hbar}}\Psi)=a^2e^{\frac{ax}{\hbar}}\Psi.
$$
More precisely, we want to look for the asymptotic of these eigenfunctions
when both $\hbar$ and $\kap$ go to $0$.
\prop{nonstat}
\begin{enumerate}
\item
There exists unique solution $\Psi$ of \refe{non-st} in
$W(\kap)$ such that $\Psi=1+O(Q)$.
\item
This solution $\Psi$ takes the form
\eq{PsiPhig}
\Psi=e^{\frac{\Phi}{\kap}+g}
\end{equation}
where $g\in Q\CC(a,\hbar)[e^x,e^{-x}][[\kap,Q]]$ and
$\Phi\in Q\CC(a,\hbar)[[\kap,Q]]$.
\item
One has
$$
\hbar Q\frac{\d \Phi}{\d Q}=b.
$$
\item
The limit
$$
\calF^{inst}=\underset{\hbar\to 0}\lim \, \hbar \Phi(a,\hbar,Q)
$$
exists and one has
\eq{renorm}
 Q\frac{\d \calF^{inst}}{\d Q}=v
\end{equation}
\end{enumerate}
\eprop

{\bf Remark.} We will explain the origin of the notation a little later.

\prf
Let us first prove (1).
Let us write
$$
\Psi=\sum\limits_{n=0}^{\infty}\Psi_n Q^n,\ \Psi_0=1.
$$
Then \refe{non-st} becomes equivalent to the sequence of equations:
\eq{Psin}
\hbar^2\Psi_n''+2\hbar a\Psi_n'-\kap n\Psi_n=U(x)\Psi_{n-1}.
\end{equation}
Let $D_n$ denote the differential operator
$\hbar^2\frac{d^2}{dx^2}+2\hbar a\frac{d}{dx}-\kap n$. Then it is easy to see that
$D_n$ is invertible when acting on $\CC(a,\hbar,\kap)[e^x,e^{-x}]$ (it is diagonal
in the basis given by the functions $\{ e^{kx}\}_{k\in\ZZ}$ with non-zero eigenvalues).
Thus by induction we get a unique solution for each $\Psi_n$, $n\geq 1$.

Let $F=\ln \Psi$. First of all, we claim that $\kap F$ is regular when $\kap =0$. This is proved exactly
in the same way as part (2) of \refp{stat} and we leave it to the reader. Let us now write
$$
F=\sum\limits_{n=-1}^{\infty} F_n \kap^n.
$$
We want to compute $F_{-1}$.

Equation \refe{non-st} is equivalent to the equation
\eq{non-st-F}
\hbar^2((F')^2+F'')+2a\hbar F'+QU(x)=\hbar\kap Q\frac{\partial F}{\partial Q}.
\end{equation}
Decomposing this in a power series in $\kap$ and looking at the coefficient of $\kap^{-2}$
we see that $\Phi=F_{-1}$ satisfies the equation
$(\Phi')^2=0$; in other words $\Phi$ is indeed independent of $x$.

Let us now look at the free term (in $\kap$) in the above identity
(it is easy to see that the coefficient of $\kap^{-1}$ is
automatically 0 on both sides). We get the equation \eq{rrr}
\hbar^2 (F_0')^2+\hbar^2 F_0''+2a\hbar F_0'+QU(x)=\hbar
Q\frac{\partial \Phi}{\partial Q}.
\end{equation}

Note now that \refe{rrr} is basically the same equation as
\refe{phi}
if we set
$b=\hbar Q\frac{\partial \Phi}{\partial Q}$ and
$F_0=\hbar^{-1}\phi$.
Since obviously $F_0|_{Q=0}=0$,
the uniqueness statement from \refp{stat}(1) implies (3). Now (4) is equivalent to \refp{stat}(3).

\epr
\defe{prepotential}
 The function
$\calF^{inst}(a,Q)$ is called {\it the instanton part} of the
prepotential.
\edefe

\medskip
\noindent
{\bf Remark.} In the context of integrable systems one is usually interested in the {\it full}
Seiberg-Witten prepotential $\calF$ which is defined as the sum of $\calF^{inst}$ and $\calF^{pert}$; here
$\calF^{pert}$ is called the perturbative part of the prepotential and it is usually given by some
simple formula. We don't know if there is a canonical choice of $\calF^{pert}$ in our generality. However,
we may observe that in all the known cases $\calF^{pert}$ satisfies the equation
$$
Q\frac{\partial \calF^{pert}}{\partial Q}=a^2.
$$
This fixes $\calF^{pert}$ uniquely up to adding a function which is independent of $Q$. Note that if we now define
$\calF=\calF^{inst}+\calF^{pert}$ (for any choice of $\calF^{pert}$ satisfying the above equation) then the equation
\refe{renorm} gets simplified: it is now equivalent to
\eq{renorm'}
Q\frac{\partial \calF}{\partial Q}=u.
\end{equation}


\sec{highdim}{Schr\"odinger operators in higher dimensions and integrable systems}

We now want to generalize the results of the previous section to
higher dimensional situation.
\ssec{}{The setup}In this section we are going to
work with the following general setup. Let $\grh$ be a finite
dimension vector space over $\CC$ and let $\Lam\subset \grh$ be a
lattice. We denote by $H$ the algebraic torus whose lattice of
co-characters is $\Lam$ (analytically one may think of $H$ as
$\grh/2\pi i\Lam$ by means of the map $x\mapsto e^{x}$); we let
$\CC[H]$ denote the algebra of polynomial functions on $H$; we might
think of elements of $\CC[H]$ as trigonometric polynomials on
$\grh$. We assume that $\grh$ is endowed with a non-degenerate
bilinear form $\la\cdot,\cdot\ra$ which takes integral values on
$\Lam$.

Let $\calK$ denote the field of rational functions on $\grh^*\x
\CC^2$ (typically, we denote an element in $\grh^*\x \CC^2$ by
$(a,\hbar,\kap)$ with $a\in\grh^*$; so, sometimes we shall write
$\CC(a,\hbar,\kap)$ instead of $\calK$).
Let $Q$ be another indeterminate. We are going to be
interested in the space $W(\kappa):=\calK[H][[Q]]$; its elements are power
series in $Q$ whose coefficients  lie in $\calK$.

Let $\Del$ denote the Laplacian on $\grh$ (or $H$) corresponding to the
bilinear form fixed above. Fix now any $P\in \CC[H]$. We shall denote by $U$ the corresponding
function on $\grh$ given by the formula
$$
U(x)=P(e^x).
$$
Now, following the previous section we define the operators
$$
T=\hbar^2\Del +QU(x);\qquad T^a=\hbar^2\Del+ 2\hbar\la\nabla,a\ra
+QU(x).
$$
Here for a function $\psi$ we denote by $\nabla \psi$ its
differential in the $\grh$-direction. Similarly, we define
$$
\calL=T-\kap Q\frac{\partial}{\partial
Q}\qquad\text{and}\qquad\calL^a=T^a-\kap Q\frac{\partial}{\partial
Q}
$$
Here $a\in \grh^*$.
Note that as before for a fixed $a$ the operator $T^a+\la a,a\ra$ is formally
conjugate to $T$, and the operator $\calL^a+\la a,a\ra$ is formally conjugate
to $\calL$.

As before, we set $W=\CC(a,\hbar)[H][[Q]]$.
Then with such notations \refp{stat} and \refp{nonstat} hold
as stated in the current
situation as well.
The proofs are just word-by-word repetitions of those from the
one-dimensional situation.

However, generalizing the results of \refss{exp1d} turns out to be
a little bit more tricky. In order to do this we need to make some
integrability assumptions.
\ssec{}{Integrability}Let us denote by $\calD$ the subalgebra of the algebra of differential operators on $H$ with
coefficients in $\CC[\hbar,Q]$ consisting of all differential operators of the form
$\sum \hbar^i D_i$ (where $D_i$ is a differential operator on $H$ with coefficients in $\CC[Q]$) such
that the order of $D_i$ is $\leq i$. It is clear that $\calD/\hbar\calD$ is canonically isomorphic
to $\calO(T^*H)\ten \CC[Q]=\calO(T^*H\x \CC)$ (here $T^*H$ denotes the cotangent bundle to $H$ and $\calO(T^*H)$ is the
algebra of regular functions on it). Note that $T^*H=H\x \grh^*$. We shall denote the resulting
map from $\calD$ to $\calO(T^*H)\ten\CC[Q]$ by $\sig$ and call it the symbol map.

Similarly we let $\calD^a=\calD\ten \calO(\grh^*)$; we have $\calD^a/\hbar\calD^a\simeq \calO(T^*H\x \CC\x\grh^*)$.
We let $\sig^a:\calD^a\to \calO(T^*H\x \CC\x\grh^*)$ denote the corresponding symbol map.

From now on we want to change our point of view a little bit and think about $T^a$ as a differential
operator on $H$ rather than on $\grh$.
Note that if we do so then $T^a$  lies in $\calD^a$.

We now assume that in addition to the above data we are given the following:

a) An affine algebraic variety $S$ such that $\dim S=\dim H$;

b) A finite morphism $\pi:\grh^*\to S$;

c) An injective homomorphism $\eta:\calO(S)\to\calD$.

\noindent
These data must satisfy the following conditions:

1) $T$ lies in the image of $\eta$; we let $C\in\calO(S)$ denote the (unique) function for which
$\eta(C)=T$.

2) $\eta|_{Q=0}$ is equal to the composition of $\pi^*:\calO(S)\to \calO(\grh^*)$ with the natural
embedding $\calO(\grh^*)\to\calD$ which sends every function $h\in\calO(\grh^*)$ which is homogeneous
of degree $d$ to $\hbar^d D_h$ where $D_h$ is the differential operator with constant coefficients corresponding
to $h$.

In this case we shall say that $T$ is {\it integrable}. Note that if
$\dim H=1$ then $T$ is automatically integrable.

Let $\eta^a:\calO(S)\to \calD^a$ denote the composition of $\eta$ with the conjugation by $e^{\frac{\la a,x\ra}{\hbar}}$.
Note that $T^a=\eta^a(C)-\la a,a\ra$.

Let also $p:T^*H\x \CC\to S$ denote the morphism such that for every $f\in \calO(S)$ we have
$$
p^*(f)=\sig\circ\eta(f).
$$
This morphism represents the classical integrable system, which is
the classical limit of the quantum integrable system defined by
$\eta$.

\ssec{}{Computation of $\lim\limits_{\hbar\to 0}
b$ via periods in the integrable case}We now want to explain how to generalize the results \refss{exp1d}
to our multi-dimensional situation in the integrable case.

First of all, the operator $T^a$ has simple spectrum in $W$; therefore the function $\psi$ which is an eigenfunction
of $T^a$, is automatically an eigenfunction of every operator of the form $\eta^a(f)$ ($f\in \calO(S)$). More precisely,
we get a homomorphism $\bfb:\calO(S)\to \calO(\grh^*)[[Q]]$ such that for each $f\in \calO(S)$ we have
$$
\eta^a(f)(\psi)(t,a,Q,\hbar)=f(\bfb(a,Q,\hbar))\psi(t,a,Q,\hbar).
\footnote{Here $t\in H$ (i.e. we think about $\psi$ as a function on $H$ rather than on $\grh$).}
$$
Note that $b=\bfb^*(C)-\la a,a\ra$. It is easy to see that the limit
$\lim\limits_{\hbar\to 0}{\bold b}^*(C)$ exists; we denote it by
$\bfu$.
 By the definition $\bfu$ is a map from $\grh^*\x\Sig$ to $S$ where $\Sig$ denotes the formal disc with coordinate
 $Q$. It is clear that $\bfu|_{Q=0}=\pi$.

Let us now look at the function ${\varphi}=\lim\limits_{\hbar\to 0} (\hbar\ln \psi)$. Then we have
\eq{fiber}
p(d{\varphi}(t,a,Q)+a,Q)=\bfu(a,Q).
\end{equation}

On the other hand, for any $\lam\in\Lam$ considered as a morphism $\lam:\CC^*\to H$ we must
have
\eq{oi}
\oint \lam^*d{\varphi}=0
\end{equation}
where $\oint$ denotes the integral over the unit circle in $\CC^*$. Let us think of $d{\varphi}$ as a morphism
$H\to T^*H$ (which depends on $a$ and $Q$). We denote by $\alp$ the canonical one-form on $T^*H$.
Let also $L_\lam$ denote the image of the unit circle under $\lam$. Then \refe{oi} is equivalent to
\eq{oi'}
\oint\limits_{L_\lam} (d{\varphi})^*\alp=0.
\end{equation}
We can now again write $a$ locally as a function of $\bfu$ and $Q$; $a=a(\bfu,Q)$. To do this we must make
a (local) choice of $a_0:=a_0(\bfu,Q)$. Note that $a_0$ must satisfy
$$
\pi(a_0)=\bfu
$$
and therefore choosing $a_0$ amounts to choosing a local branch of $\pi$. Let now $\lam$ be as above. Then we denote
by $A_{\lam,\bfu,Q}$ the unique 1-dimensional cycle in $T^*H$ such that:

1) The projection of $A_{\lam,\bfu,Q}$ to $H$ is equal to $L_\lam$;

2) $A_{\lam,\bfu,Q}\subset p^{-1}(\bfu)$;

3) $A_{\lam,\bfu,Q}$ depends continuously on $Q$ and for $Q=0$ it
lies in the above chosen branch of $\pi$. \footnote{More
precisely, this means the following: for $Q=0$ the map $p$ is
equal to the composition of the natural projection $T^*H\to
\grh^*$ and $\pi:\grh^*\to S$. Thus for every $\bfu$ we have
$p|_{Q=0}^{-1}(\bfu)=H\x \pi^{-1}(\bfu)$. We require that $A_\lam$
lie in the product of $H$ and the corresponding branch of $\pi$.}

Then \refe{fiber} says that for every $\lam\in\Lam$ we have
\eq{aeq} \la a,\lam\ra=\frac{1}{2\pi
i}\oint\limits_{A_{\lam,\bfu,Q}}\alp.
\end{equation}

\ssec{}{Some variants} Let us choose a closed cone
$\grh^*_+\subset \grh^*$ which is integral with respect to $\Lam$
(i.e. given by finitely many inequalities given by elements of
$\Lam$). We assume also that $a\in \grh^*_+$ implies that
$-a\not\in\grh_+^*$ for $a\neq 0$ (i.e. $0$ is an extremal point
of $\grh_+$). Set $\Lam^{\vee}_+=\Lam^{\vee}\cap \grh_+^*$. We
denote by $\hatW$ the corresponding completion of $W$; by the
definition it consists of all formal sums
$$
\sum c_\gam e^{\la \gam,x\ra}
$$
where $\gam\in\Lam^{\vee}$ and such that for each $\chlam\in\Lam^{\vee}$ the set
$$
\{ \gam\in \chlam-\Lam_+|\ \text{such that } c_{\gam}\neq 0\}
$$
is finite.

It is easy to see that the results of this section generalize immediately
to the situation when the initial Schr\"odinger operator $T$ takes the form
$$
T=\hbar^2\Del +\bfU(Q,x)
$$
where $\bfU\in \CC[H][[Q]]$ subject to the following condition:

\medskip

$\bullet$ The function
$\bfU(0,x)$ is a linear combination of $e^{\la \chlam,x\ra}$ with $\chlam\in \grh^*_+,\ \chlam\neq 0$.

\medskip

In this case the eigenfunctions $\psi$ and $\Psi$ should be elements of respectively $\hatW$ and $\hatW(\kap)$.

The above condition guarantees
in particular that $0$ is an eigenvalue of $T^a$ on $\hatW$. The definition of the prepotential goes through in this case
without any changes.

Here is the basic example of the above situation. Let
$$
\grh=\{ (x_1,...,x_n)\in \CC^n\}/\CC(1,...,1)\, \quad \Lam=\{ (x_1,...,x_n)\in \ZZ^n\}/\ZZ(1,...,1).
$$
Clearly,
$$
\grh^*=\{ (a_1,...,a_n)|\ \sum a_i=0\}
$$
and we set
$$
\grh^*_+=\{(a_1,...,a_n)\in\grh^*|\ a_1+a_2+...+a_k\geq 0\ \text{for each }1\leq k\leq n\}.
$$
Let
$$
\bfU(Q,x)=2(e^{x_1-x_2}+e^{x_2-x_3}+...+e^{x_{n-1}-x_n}+Qe^{x_n-x_1})
$$
be the periodic Toda potential. It is clear that the condition above is satisfied and therefore we may speak of the
corresponding prepotential. In the next section we explain its connection with the standard physical
definition of the prepotential.

The periodic Toda potential is equal to the Toda potential defined by
the affine Lie algebra $\widehat {sl_n}$ (cf. for example
\cite{etingof}). One can easily see that the Toda potential for any affine Lie algebra (cf. \cite{etingof}) satisfies our
conditions and thus the corresponding prepotential is well-defined.

Note also that the operator
$$
\hbar^2\Del+\bfU(Q,x)
$$
turns into the operator
$$
\hbar^2\Del+2Q^{1/n}(e^{x_1-x_2}+...+e^{x_{n-1}-x_n}+e^{x_n-x_1}).
$$
after the change of variables
$$
x_j\mapsto x_j+\frac{j\ln Q}{n}.
$$

Thus when computing the prepotential we may deal with the latter operator
(a similar statement is true for all affine Lie algebras).

\medskip
\noindent
{\bf Remark.} The variable $Q$ that we are using is connected with the variable $\Lam$ (which is commonly
used by physicists - cf. \cite{NaYo1}, \cite{nek} etc.) by the formula
$$
Q=\Lam^{2n}.
$$

\medskip
\noindent {\bf Remark.} It is not difficult to check that  if $U$
is equal to the Toda potential for $\widehat{sl_n}$ then our
definition of $\calF^{inst}$ coincides with the one usually given
by physicists (cf. Chapter 2 of \cite{NaYo1}). Let us give a very
brief sketch of the proof of this result (details will appear in a
subsequent publication in a more general setting). Namely, in this
case \refe{renorm} becomes equivalent to the {\it renormalisation
group equation} (Proposition 2.10 of \cite{NaYo1}). Note that in
the original (Seiberg-Witten) definition of the prepotential $a$
is expressed in terms of periods of some family of curves over
$\grh^*/W\x \CC$ where $W=S_n$ is the Weyl group of $sl_n$ (here
the second factor is the line with coordinate $Q$). These curves
are called the {\it Seiberg-Witten curves} (cf. Section 2.1 of
\cite{NaYo1}). In our case, $a$ is expressed via periods on the
fibers of the map $p_Q:T^*H\to \grh^*/W$ (note that in this case
we have $S=\grh^*/W$). However, it is well known (cf.
\cite{krich}) that the fibers of the map $p_Q$ are open pieces in
the Jacobians of the Seiberg-Witten curves; thus periods of a
regular one-form over these fibers are equal to the periods of a
certain meromorphic one-form over the curves themselves and it is
not difficult to check that we get exactly the same periods as we
need. Some generalization of this fact will be considered in much
more detail in a further publication.

\sec{}{Proof of Nekrasov's conjecture}

In this section we want to prove \reft{main'} (and thus also \reft{main}).
The first part of \reft{main'} is an immediate corollary of Corollary 3.7 from
\cite{Br} combined with the definition of $\calF^{inst}$ given by \refd{prepotential}.
Thus it remains to prove the second part of \reft{main'}. The proof is based on
the following result.

\th{pg} Let $P\subset G$ be a parabolic and
let $(d,\theta)\in \Lam_{G,P}^{\aff,+}$. Then one of the following is true:

a) Both $\int\limits_{\calU_{G,P}^{d,\theta}}1$and ${\int\limits_{\calU^d_G}1}$ are 0.

b) ${\int\limits_{\calU^d_G}1}\neq 0$ and the ratio
$$
\frac{\int\limits_{\calU_{G,P}^{d,\theta}}1}{\int\limits_{\calU^d_G}1}
$$
is regular when $\eps_2\to 0$.
\eth

Let us first explain why \reft{pg} implies \reft{main'}.
First of all, we claim that
for any $d\leq d'$ we have
$$
{\int\limits_{\calU^{d'}_G}1}=A_d {\int\limits_{\calU^d_G}1}
$$
where $A_d$ is a regular function on $\grh\times \CC^2$ (in particular, it is regular
when $\eps_2\to 0$). Indeed, according  to \cite{bfg} there exists a closed
$G\x(\CC^*)^2$-equivariant embedding
$\calU^d\overset{i_d}\to \calU^{d'}$. Since $\calU_{d'}$ is contractible we have
$H^*_{G\x(\CC^*)^2}(\calU_G^{d'})=\calA_{G\x (\CC^*)^2}$. Thus it follows that
the direct image $(i_d)_*1$ of the equivariant unit cohomology class is equal to
some $A_d\in \calA_{G\x (\CC^*)^2}$.

Now it follows
from \reft{pg} that the ratio
$$
\frac{\calZ_{G,P}^{\aff}(\grq,Q,a,\eps_1,\eps_2)}{\calZ_G^{\aff}(Q,a,\eps_1,\eps_2)}
$$
is regular when $\eps_2\to 0$. This means that
$$
\underset{\eps_2\to 0}\lim \eps_2(\ln \calZ_{G,P}^{\aff} -
\ln\calZ_G^{\aff})=0.
$$
This is the statement of \reft{main'}(2).

Thus to complete the proof we need to prove \reft{pg}. The proof  is
based on the following general lemma.
\lem{general}
Let $L_1$ and $L_2$ be two algebraic tori and let $L=L_1\x L_2$.
We let $\grl_1$ and $\grl_2$ denote the corresponding Lie
algebras.
We shall denote a typical element in $\grl$ by $(l_1,l_2)$, where
$l_i\in\grl_i$.

Let  $\pi:X\to Y$ be a morphism of $L$-varieties. Assume that:

1) Both $X^L$ and $Y^L$ are proper.

2) The natural map $X^{L_1} \to Y^{L_1}$ is proper.

Then if ${\int\limits_Y 1}$ is zero then ${\int\limits_X 1}$ is also zero (here we consider both integrals in
$L$-equivariant cohomology). If ${\int\limits_Y 1}\neq 0$ then
the ratio
$$
\frac{\int\limits_X 1}{\int\limits_Y 1}
$$
(where the integral is taken in $L$-equivariant cohomology) is
regular when $l_2\to 0$.
\elem
\refl{general} is an easy corollary of the definition of the above integrals given in Section 2 of
\cite{Br} and we leave the proof to the reader.
\ssec{}{End of the proof}First of all, we need to show that
$$
\int\limits_{\calU^d_G}1\neq 0.
$$

We now want to apply \refl{general} to the
case when $X=\calU_{G,P}^{d,\theta}$, $Y=\calU_G^d$, $L_1=T\x \CC^*$
where
the $\CC^*$ factor corresponds to $\eps_1$ and $L_2=\CC^*$
corresponding to $\eps_2$. To avoid confusion in the notation we shall
denote
the ``first'' (i.e. horizontal) copy  of $\CC^*$ by $\CC^*_1$ and the
other copy by $\CC^*_2$. We need to show that
the map $(\calU_{G,P}^{d,\theta})^{T\x \CC^*_1}\to (\calU_G^d)^{T\x
\CC^*_1}$
is proper. In fact, we claim that the following stronger statement is
true:
\lem{}
The map $(\calU_{G,P}^{d,\theta})^{\CC^*_1}\to (\calU_G^d)^{\CC^*_1}$
is an isomorphism.
\elem
\prf
This is an easy corollary of Theorem 10.2 of \cite{bfg}. In {\it loc. cit} a
natural stratification of $\calU_{G,P}$ is described and it follows
immediately that
$$
(\calU_{G,P}^{d,\theta})^{\CC^*_1}=
(\calU_G^d)^{\CC^*_1}=\Sym^d(\bfX\backslash\{ \infty\})
$$
(recall that $\bfX$ denotes the ``vertical'' axis in $\PP^2$).
\epr


\end{document}